\documentclass[a4paper,twoside,11pt]{article}
\usepackage{amsfonts}
\usepackage{amssymb}\usepackage{amsmath}
\newtheorem{thm}{Theorem}[section]
 \newtheorem{cor}{Corollary}[section]
 \newtheorem{lem}[thm]{Lemma}

 \numberwithin{equation}{section}
\newcommand{\double}{\baselineskip 1.24 \baselineskip}

\title{On  an improvement of the Hardy-Hilbert type  inequality
}

\author{{Guang-Sheng  Chen\thanks{\text{E-mail address}: cgswavelets@126.com(Chen)
}\quad}\\
{\small Department of Computer Engineering, Guangxi Modern
Vocational Technology College,} \\{\small Hechi,Guangxi, 547000,
P.R. China}
}
\begin{document}
\date{}
\maketitle \double

\textbf{Abstract:}\quad In this paper, by estimating
 the weight coefficient effectively, we establish an
improvement of a Hardy-Hilbert type inequality proved by B.C.
Yang, our main tool is Euler-Maclaurin expansion for the zeta function. As applications, some particular results are considered.   

\section{Introduction}
\hskip\parindent

 If $p,q >1$, $\textstyle{1 \over p} + \textstyle{1 \over q} = 1$, $a_n $,
$b_n \ge 0$, $0 < \sum\limits_{n = 1}^\infty {a_n^p } < \infty
$,and $0 < \sum\limits_{n = 1}^\infty {b_n^q } < \infty $, then[1]
have
\begin{equation*}
\sum\limits_{n = 1}^\infty  {\sum\limits_{m = 1}^\infty  {\frac{{{a_m}{b_n}}}
{{m + n}}} }  < \frac{\pi }
{{\sin (\pi /p)}}{\left\{ {\sum\limits_{n = 1}^\infty  {a_n^p} } \right\}^{\frac{1}
{p}}}{\left\{ {\sum\limits_{n = 1}^\infty  {b_n^q} } \right\}^{\frac{1}
{q}}}
, \label{1.1}\tag{1.1}
\end{equation*}
\begin{equation*}
 \sum\limits_{n = 1}^\infty {\sum\limits_{m = 1}^\infty
{\frac{a_m b_n }{\max \{m,n\}}} } < pq\left\{\sum\limits_{n = 1}^\infty
{a_n^p } \right\}^{\frac{1}{p}}\left\{\sum\limits_{n = 1}^\infty {b_n^q }
\right\}^{\frac{1}{q}}, \label{1.2}\tag{1.2}
\end{equation*}
where the constant factor $\frac{\pi }{\sin (\pi / p)}$ and $pq$
are the best possible. Inequality (1.1) is well known as Hardy-
Hilbert's inequality and (1.2) is named a Hardy- Hilbert's type
inequality. Both of them are important in analysis and its
applications [2]. In the recent years, a lot of results with
generalizations of this type of inequality were obtained (see
[3]). Under the same conditions as (1.1) and (1.2), there are some
Hardy-Hilbert's type inequalities similar to (1.1) and (1.2),
which also had been studied and generalized by some
mathematicians.

Recently, by introducing a parameter Yang [4] gave the following
generalization of inequality (1.2) for $p=q=2$:

If $p,q >1$, $\textstyle{1 \over p} + \textstyle{1 \over q} = 1$,
$2 - \min \{p,q\} < \lambda \le 2$, $a_n $, $b_n \ge 0$, such that
$0 < \sum\limits_{n = 1}^\infty {n^{(p - 1)(2 - \lambda ) -
1}a_n^p } < \infty $, and $0 < \sum\limits_{n = 1}^\infty {n^{(q -
1)(2 - \lambda ) - 1}b_n^q } < \infty $, then
\begin{equation*}
\sum\limits_{n = 1}^\infty {\sum\limits_{m = 1}^\infty {\frac{a_m
b_n }{\max \{m^\lambda ,n^\lambda \}}} } < k_\lambda (p)\left\{
{\sum\limits_{n = 1}^\infty {n^{(p - 1)(2 - \lambda ) - 1}a_n^p }
} \right\}^{\textstyle{1 \over p}}\left\{ {\sum\limits_{n =
1}^\infty {n^{(q - 1)(2 - \lambda ) - 1}b_n^q } }
\right\}^{\textstyle{1 \over q}},\label{1.3}\tag{1.3}
\end{equation*}
where the constant factor $k_\lambda (p) = \frac{\lambda pq}{(p +
\lambda - 2)(q + \lambda - 2)}$ is the best possible.

In this paper, by introducing a parameter and estimating the
weight coefficient, we obtain an improvement of (1.3).As
application, we consider some particular results.

\section{Preliminary results }
\hskip\parindent
First, we need the following formula of the Riemann-$\zeta $
function (see [5]):
\begin{equation*}
 \zeta (\rho ) = \sum\limits_{n = 1}^m {\frac{1}{n^\rho }} - \frac{m^{1 -
\rho }}{1 - \rho } - \frac{1}{2m^\rho } - \sum\limits_{n = 1}^{l -
1} {\frac{B_{2n} }{2n}\left( {{\begin{array}{*{20}c}
 { - \rho } \hfill \\
 {2n - 1} \hfill \\
\end{array} }} \right)\frac{1}{m^{\rho + 2n - 1}}}
- \frac{B_{2l} }{2l}\left( {{\begin{array}{*{20}c}
 { - \rho } \hfill \\
 {2l - 1} \hfill \\
\end{array} }} \right)\frac{\varepsilon }{m^{\rho + 2l - 1}}, \\
 \label{2.1}\tag{2.1}
\end{equation*}
 where $\rho > 0$, $\rho \ne 1$, $m$, $l \ge 1$,
$m,l \in {\rm N}$, $0 < \varepsilon = \varepsilon (\rho ,l,m) <
1$.The numbers $B_1 = - 1 / 2$, $B_2 = 1 / 6$, $B_3 = 0$, $B_4 = -
1 / 30$，\ldots are Bernoulli numbers. In particular, $\zeta (\rho
) = \sum\limits_{n = 1}^\infty {\frac{1}{n^\rho }} $ ($\rho > 1$).

Since $\zeta (0) = - 1 / 2$, then the formula of the
Riemann-$\zeta $ function (2.1) is also true for $\rho = 0$.
\begin{lem}\label{Lemma 2.1.}
If $p,q >1$, $\textstyle{1 \over p} + \textstyle{1 \over q} = 1$,
$2 - \min \{p,q\} < \lambda \le 2$, define the weight coefficients
$\omega (m,\lambda ,p)$ and $\omega (n,\lambda ,q)$ as
\begin{equation*}
\omega (m,\lambda ,p) = \sum\limits_{n = 1}^\infty  {\frac{1}
{{\max \{ {m^\lambda },{n^\lambda }\} }}{{\left( {\frac{m}
{n}} \right)}^{(2 - \lambda )/p}}} ,
\label{2.2}\tag{2.2}
\end{equation*}
\begin{equation*}
\omega (n,\lambda ,q) = \sum\limits_{m = 1}^\infty  {\frac{1}
{{\max \{ {m^\lambda },{n^\lambda }\} }}{{\left( {\frac{n}
{m}} \right)}^{(2 - \lambda )/q}}} ,
\label{2.3}\tag{2.3}
\end{equation*}
then we have
\begin{equation*}
\omega (m,\lambda ,p) < {m^{1 - \lambda }}\left[ {{k_\lambda }(p) - \frac{p}
{{3(p + \lambda  - 2){m^{(p + \lambda  - 2)/p}}}}} \right],
\label{2.4}\tag{2.4}
\end{equation*}
and
\begin{equation*}
\omega (n,\lambda ,q) < {n^{1 - \lambda }}\left[ {{k_\lambda } - \frac{q}
{{3(p + \lambda  - 2){n^{(q + \lambda  - 2)/q}}}}} \right],
\label{2.5}\tag{2.5}
\end{equation*}
where $k_\lambda = \frac{\lambda pq}{(p + \lambda - 2)(q + \lambda
- 2)}$.
\end{lem}
\textbf{Proof.} For $2 - \min \{p,q\}< \lambda \le 2$, taking
$\rho = \textstyle{{2 - \lambda } \over p} \ge 0$, $l = 1$ in
(2.1),we get
\begin{equation*}
\zeta (\tfrac{{2 - \lambda }}
{p}) = \sum\limits_{n = 1}^m {\frac{1}
{{{n^{(2 - \lambda )/p}}}}}  - \frac{{p{m^{(p + \lambda  - 2)/p}}}}
{{p + \lambda  - 2}} - \frac{1}
{{2{n^{(2 - \lambda )/p}}}} + \frac{{2 - \lambda }}
{{12p{m^{1 + (2 - \lambda )/p}}}}{\varepsilon _1} ,\label{2.6}\tag{2.6}
\end{equation*}
Where $0 < \varepsilon _1 < 1$.

Taking $\rho = \textstyle{2 \over p} + \textstyle{\lambda \over
q}$,we obtain
\begin{equation*}
\zeta (\tfrac{2}
{p} + \tfrac{\lambda }
{q}) = \sum\limits_{n = 1}^{m - 1} {\frac{1}
{{{n^{2/p + \lambda /q}}}}}  + \frac{{q{m^{ - (q + \lambda  - 2)/q}}}}
{{q + \lambda  - 2}} + \frac{1}
{{2{m^{2/p + \lambda /q}}}} + \frac{{p\lambda  + 2q}}
{{12pq{m^{1 + 2/p + \lambda /q}}}}{\varepsilon _2} , \label{2.7}\tag{2.7}
\end{equation*}
Where $0 < \varepsilon _2 < 1$.

Thus we get
\begin{equation*}
\begin{split}
& \omega (m,\lambda ,p) = \sum\limits_{n = 1}^\infty  {\frac{1}
{{\max \{ {m^\lambda },{n^\lambda }\} }}{{\left( {\frac{m}
{n}} \right)}^{(2 - \lambda )/p}}} \\
  & = \sum\limits_{n = 1}^m {\frac{1}
{{\max \{ {m^\lambda },{n^\lambda }\} }}{{\left( {\frac{m}
{n}} \right)}^{(2 - \lambda )/p}}}  - \frac{1}
{{{m^\lambda }}} + \sum\limits_{n = m}^\infty  {\frac{1}
{{\max \{ {m^\lambda },{n^\lambda }\} }}{{\left( {\frac{m}
{n}} \right)}^{(2 - \lambda )/p}}} \\
   &= \sum\limits_{n = 1}^m {\frac{1}
{{{m^\lambda }}}{{\left( {\frac{m}
{n}} \right)}^{(2 - \lambda )/p}}}  - \frac{1}
{{{m^\lambda }}} + \sum\limits_{n = m}^\infty  {\frac{1}
{{{n^\lambda }}}{{\left( {\frac{m}
{n}} \right)}^{(2 - \lambda )/p}}} \\
  & = \frac{1}
{{{m^{(p + 1)\lambda  - 2}}}}\sum\limits_{n = 1}^m {\frac{1}
{{{n^{(2 - \lambda )/p}}}}}  - \frac{1}
{{{m^\lambda }}} + {m^{(2 - \lambda )/p}}\sum\limits_{n = m}^\infty  {\frac{1}
{{{n^{2/p + \lambda /q}}}}}. \\
\end{split} \notag
\end{equation*}
By (2.6) and (2.7), we have
\begin{equation*}
\begin{split}
 & \omega (m,\lambda ,p) < \frac{1}
{{{m^{(p + 1)\lambda  - 2}}}}\left[ {\zeta (\tfrac{{2 - \lambda }}
{p}) + \frac{{p{m^{(p + \lambda  - 2)/p}}}}
{{p + \lambda  - 2}} + \frac{1}
{{2{m^{(2 - \lambda )/p}}}}} \right] - \frac{1}
{{{m^\lambda }}} \hfill \\
  & + {m^{(2 - \lambda )/p}}\left[ {\frac{{q{m^{ - (q + \lambda  - 2)/q}}}}
{{q + \lambda  - 2}} + \frac{1}
{{2{m^{2/p + \lambda /q}}}} + \frac{{p\lambda  + 2q}}
{{12pq{m^{1 + 2/p + \lambda /q}}}}} \right] \hfill \\
  & = \frac{1}
{{{m^{(p + 1)\lambda  - 2}}}}\zeta (\tfrac{{2 - \lambda }}
{p}) + \frac{{p{m^{1 - \lambda }}}}
{{p + \lambda  - 2}} + \frac{1}
{{2{m^\lambda }}} - \frac{1}
{{{m^\lambda }}} + \frac{{q{m^{1 - \lambda }}}}
{{q + \lambda  - 2}} + \frac{1}
{{2{m^\lambda }}} + \frac{{p\lambda  + 2q}}
{{12pq{m^{1 + \lambda }}}} \hfill \\
   &= \frac{1}
{{{m^{(p + 1)\lambda  - 2}}}}\zeta (\tfrac{{2 - \lambda }}
{p}) + \frac{{\lambda pq{m^{1 - \lambda }}}}
{{(p + \lambda  - 2)(q + \lambda  - 2)}} + \frac{{p\lambda  + 2q}}
{{12pq{m^{1 + \lambda }}}} \hfill \\
  & = {m^{1 - \lambda }}\left\{ {\frac{{\lambda pq}}
{{(p + \lambda  - 2)(q + \lambda  - 2)}} - \frac{1}
{{{m^{(p + \lambda  - 2)/p}}}}\left[ { - \zeta (\tfrac{{2 - \lambda }}
{p}) - \frac{{p\lambda  + 2q}}
{{12pq{m^{(p + \lambda  - 2)/p}}}}} \right]} \right\}.\\
\end{split} \notag
\end{equation*}
In (2.6), taking $m = 1$,by $2 - \min \{p,q\} < \lambda \le 2$, we
obtain
\begin{equation*}
\begin{split}
 &\zeta (\textstyle{{2 - \lambda } \over p}) = 1 - \frac{p}{p + \lambda - 2}
- \frac{1}{2} + \frac{(2 - \lambda )\varepsilon _1 }{12p} \\
 &< \frac{1}{2} - \frac{p}{p + \lambda - 2} + \frac{2 - \lambda }{12p} \\
& = - \frac{(\lambda - 2 - 3p)(\lambda - 2 - 2p)}{12p(p + \lambda
- 2)} < 0.
\end{split} \notag
\end{equation*}
Therefore for $m \ge 1$, $m \in {\rm N}$,$2 - \min \{p,q\} <
\lambda \le 2$,we have
\begin{equation*}
\begin{split}
 &- \zeta (\textstyle{{2 - \lambda } \over p}) - \frac{p\lambda +
2q}{12pqm^{\textstyle{{p + \lambda - 2} \over p}}} >
\frac{(\lambda - 2 -
3p)(\lambda - 2 - 2p)}{12p(p + \lambda - 2)} - \frac{p\lambda + 2q}{12pq} \\
 &= \frac{q(\lambda - 2 - 3p)(\lambda - 2 - 2p) - (p\lambda + 2q)(p + \lambda
- 2)}{12pq(p + \lambda - 2)} \\
& = \frac{q(\lambda - 2) - 3p) + (p\lambda + 5pq + 2q)(2 - \lambda
) -
p(p\lambda + 2q) + 6p^2q}{12pq(p + \lambda - 2)} \\
 &\ge \frac{ - (2p + 2q) + 6pq}{12pq(p + \lambda - 2)} = \frac{p}{3(p +
\lambda - 2)}. \\
\end{split} \notag
\end{equation*}
Applying the last result and the inequality for $\omega (m,\lambda
,p)$ above, we obtain (2.4).Similarly, we can prove (2.5).The
lemma is proved.

\section{Main results}

\begin{thm}
\label{Theorem 3.1} If $p,q >1$, $\textstyle{1 \over p} +
\textstyle{1 \over q} = 1$, $2 - \min \{p,q\} < \lambda \le
2$,$a_n \ge 0$, $b_n \ge 0$, such that $0 < \sum\limits_{n =
1}^\infty {n^{(p - 1)(2 - \lambda ) - 1}a_n^p } < \infty $, and $0
< \sum\limits_{n = 1}^\infty {n^{(q - 1)(2 - \lambda ) - 1}b_n^q }
< \infty $, then
\begin{equation*}
\begin{split}
&\sum\limits_{n = 1}^\infty  {\sum\limits_{m = 1}^\infty  {\frac{{{a_m}{b_n}}}
{{\max \{ {m^\lambda },{n^\lambda }\} }}} }  < {\left\{ {\sum\limits_{n = 1}^\infty  {\left[ {{k_\lambda } - \frac{p}
{{3(p + \lambda  - 2){n^{(p + \lambda  - 2)/p}}}}} \right]{n^{(p - 1)(2 - \lambda ) - 1}}a_n^p} } \right\}^{\tfrac{1}
{p}}} \hfill \\
 & \quad \quad \quad \quad \quad \quad \quad \quad  \times {\left\{ {\sum\limits_{n = 1}^\infty  {\left[ {{k_\lambda } - \frac{q}
{{3(q + \lambda  - 2){n^{(q + \lambda  - 2)/q}}}}} \right]{n^{(q - 1)(2 - \lambda ) - 1}}b_n^q} } \right\}^{\tfrac{1}
{q}}}, \\
\end{split}
\label{3.1}\tag{3.1}
\end{equation*}
\begin{equation*} {\sum\limits_{n = 1}^\infty  {{n^{p + \lambda  - 3}}\left[ {\sum\limits_{m = 1}^\infty  {\frac{{{a_m}}}
{{\max \{ {m^\lambda },{n^\lambda }\} }}} } \right]} ^p} < k_{_\lambda }^{p - 1}\sum\limits_{n = 1}^\infty  {\left[ {{k_\lambda } - \frac{p}
{{3(p + \lambda  - 2){n^{(p + \lambda  - 2)/p}}}}} \right]{n^{(p - 1)(2 - \lambda ) - 1}}a_n^p}, \label{3.2}\tag{3.2}
\end{equation*}
where $k_\lambda = \frac{\lambda pq}{(p + \lambda - 2)(q + \lambda
- 2)} > 0$.
\end{thm}
\textbf{Proof.} By H\"{o}lder inequality(see[6]), we have
\begin{equation*}
\begin{split}
 &\sum\limits_{n = 1}^\infty  {\sum\limits_{m = 1}^\infty  {\frac{{{a_m}{b_n}}}
{{\max \{ {m^\lambda },{n^\lambda }\} }}} }  = \sum\limits_{n = 1}^\infty  {\sum\limits_{m = 1}^\infty  {\left[ {\frac{{{a_m}}}
{{{{(\max \{ {m^\lambda },{n^\lambda }\} )}^{1/p}}}}\frac{{{m^{(2 - \lambda )/{q^2}}}}}
{{{n^{(2 - \lambda )/{p^2}}}}}} \right]\left[ {\frac{{{b_n}}}
{{{{(\max \{ {m^\lambda },{n^\lambda }\} )}^{1/q}}}}\frac{{{n^{(2 - \lambda )/{p^2}}}}}
{{{m^{(2 - \lambda )/{q^2}}}}}} \right]} }  \hfill \\
  & \leqslant {\left\{ {\sum\limits_{n = 1}^\infty  {\sum\limits_{m = 1}^\infty  {\left[ {\frac{{a_m^p}}
{{\max \{ {m^\lambda },{n^\lambda }\} }}\frac{{{m^{p(2 - \lambda )/{q^2}}}}}
{{{n^{(2 - \lambda )/p}}}}} \right]} } } \right\}^{\tfrac{1}
{p}}}{\left\{ {\sum\limits_{n = 1}^\infty  {\sum\limits_{m = 1}^\infty  {\left[ {\frac{{{b_n}}}
{{\max \{ {m^\lambda },{n^\lambda }\} }}\frac{{{n^{q(2 - \lambda )/{p^2}}}}}
{{{m^{(2 - \lambda )/q}}}}} \right]} } } \right\}^{\tfrac{1}
{q}}} \\
  & = {\left\{ {\sum\limits_{n = 1}^\infty  {\sum\limits_{m = 1}^\infty  {\frac{1}
{{\max \{ {m^\lambda },{n^\lambda }\} }}{{\left( {\frac{m}
{n}} \right)}^{(2 - \lambda )/p}}{m^{(2 - \lambda )(p - 2)}}a_m^p} } } \right\}^{\tfrac{1}
{p}}}{\left\{ {\sum\limits_{n = 1}^\infty  {\sum\limits_{m = 1}^\infty  {\frac{1}
{{\max \{ {m^\lambda },{n^\lambda }\} }}{{\left( {\frac{n}
{m}} \right)}^{(2 - \lambda )/q}}} {n^{(2 - \lambda )(q - 2)}}b_n^q} } \right\}^{\tfrac{1}
{q}}}. \\
\end{split}
\notag
 \end{equation*}
Hence, By (2.4), (2.5), inequality (3.1) holds.

By H\"{o}lder inequality and Lemma 2.1, we obtain
\begin{equation*}
\begin{split}
 &\sum\limits_{m = 1}^\infty  {\frac{{{a_m}}}
{{\max \{ {m^\lambda },{n^\lambda }\} }}}  = \sum\limits_{m = 1}^\infty  {\frac{1}
{{{{(\max \{ {m^\lambda },{n^\lambda }\} )}^{1/p}}}}\frac{{{m^{(2 - \lambda )/{q^2}}}}}
{{{n^{(2 - \lambda )/{p^2}}}}}{a_m}\frac{1}
{{{{(\max \{ {m^\lambda },{n^\lambda }\} )}^{1/q}}}}\frac{{{n^{(2 - \lambda )/{p^2}}}}}
{{{m^{(2 - \lambda )/{q^2}}}}}} \\
  & \leqslant {\{ \sum\limits_{m = 1}^\infty  {\frac{1}
{{\max \{ {m^\lambda },{n^\lambda }\} }}\frac{{{m^{p(2 - \lambda )/{q^2}}}}}
{{{n^{(2 - \lambda )/p}}}}a_m^p} \} ^{\tfrac{1}
{p}}}{\{ \frac{1}
{{\max \{ {m^\lambda },{n^\lambda }\} }}\frac{{{n^{q(2 - \lambda )/{p^2}}}}}
{{{m^{(2 - \lambda )/q}}}}\} ^{\tfrac{1}
{q}}}\\
  & = {\{ \sum\limits_{m = 1}^\infty  {\frac{1}
{{\max \{ {m^\lambda },{n^\lambda }\} }}{{\left( {\frac{m}
{n}} \right)}^{(2 - \lambda )/p}}{m^{(2 - \lambda )(p - 2)}}a_m^p} \} ^{\tfrac{1}
{p}}}{\{ \sum\limits_{m = 1}^\infty  {\frac{1}
{{\max \{ {m^\lambda },{n^\lambda }\} }}{{\left( {\frac{n}
{m}} \right)}^{(2 - \lambda )/q}}{n^{(2 - \lambda )(q - 2)}}} \} ^{\tfrac{1}
{q}}}  \\
   &= {\{ \sum\limits_{m = 1}^\infty  {\frac{1}
{{\max \{ {m^\lambda },{n^\lambda }\} }}{{\left( {\frac{m}
{n}} \right)}^{(2 - \lambda )/p}}]{m^{(2 - \lambda )(p - 2)}}a_m^p} \} ^{\tfrac{1}
{p}}}{\{ \omega (n,\lambda ,q){n^{(2 - \lambda )(q - 2)}}\} ^{\tfrac{1}
{q}}}  \\
  & < {\{ \sum\limits_{m = 1}^\infty  {\frac{1}
{{\max \{ {m^\lambda },{n^\lambda }\} }}{{\left( {\frac{m}
{n}} \right)}^{(2 - \lambda )/p}}{m^{(2 - \lambda )(p - 2)}}a_m^p} \} ^{\tfrac{1}
{p}}}{\{ {k_\lambda }{n^{(2 - \lambda )(q - 2) + 1 - \lambda }}\} ^{\tfrac{1}
{q}}} . \\
\end{split}
\notag
 \end{equation*}
Then
\begin{equation*}
\begin{split}
& {\sum\limits_{n = 1}^\infty  {{n^{p + \lambda  - 3}}\left[ {\sum\limits_{m = 1}^\infty  {\frac{{{a_m}}}
{{\max \{ {m^\lambda },{n^\lambda }\} }}} } \right]} ^p} < k_{_\lambda }^{p - 1}\sum\limits_{n = 1}^\infty  {\sum\limits_{m = 1}^\infty  {\frac{1}
{{\max \{ {m^\lambda },{n^\lambda }\} }}{{\left( {\frac{m}
{n}} \right)}^{(2 - \lambda )/p}}{m^{(2 - \lambda )(p - 2)}}a_m^p} } \\
  & < k_{_\lambda }^{p - 1}\sum\limits_{m = 1}^\infty  {{m^{(2 - \lambda )(p - 2)}}a_m^p\sum\limits_{n = 1}^\infty  {\frac{1}
{{\max \{ {m^\lambda },{n^\lambda }\} }}{{\left( {\frac{m}
{n}} \right)}^{(2 - \lambda )/p}}} }  = k_{_\lambda }^{p - 1}\sum\limits_{m = 1}^\infty  {{m^{(2 - \lambda )(p - 2)}}\omega (m,\lambda ,p)a_m^p}. \\
\end{split}
\notag
 \end{equation*}
By Lemma 2.1, the proof of Theorem 3.1 is completed.

By Theorem 3.1, we have
\begin{cor}
\label{Corollary 3.1}   If $p,q>1$, $\textstyle{1 \over p} +
\textstyle{1 \over q} = 1$, $2 - \min \{p,q\} < \lambda \le
2$,$a_n \ge 0$, $b_n \ge 0$, such that $0 < \sum\limits_{n =
1}^\infty {n^{(p - 1)(2 - \lambda ) - 1}a_n^p } < \infty $, and $0
< \sum\limits_{n = 1}^\infty {n^{(q - 1)(2 - \lambda ) - 1}b_n^q }
< \infty $, then
\begin{equation*}
\begin{split}
 &\sum\limits_{n = 1}^\infty  {\sum\limits_{m = 1}^\infty  {\frac{{{a_m}{b_n}}}
{{\max \{ {m^\lambda },{n^\lambda }\} }}} }  < {\left\{ {\sum\limits_{n = 1}^\infty  {\left[ {{k_\lambda } - \frac{1}
{{3(p + \lambda  - 2){n^{(p + \lambda  - 2)/p}}}}} \right]{n^{(p - 1)(2 - \lambda ) - 1}}a_n^p} } \right\}^{\tfrac{1}
{p}}} \\
  &\quad \quad \quad \quad \quad \quad \quad \quad  \times {\left\{ {\sum\limits_{n = 1}^\infty  {\left[ {{k_\lambda } - \frac{1}
{{3(q + \lambda  - 2){n^{(q + \lambda  - 2)/q}}}}} \right]{n^{(q - 1)(2 - \lambda ) - 1}}b_n^q} } \right\}^{\tfrac{1}
{q}}} , \\
 \end{split}
 \label{3.3}
 \tag{3.3}
\end{equation*}
\begin{equation*}
 {\sum\limits_{n = 1}^\infty  {{n^{p + \lambda  - 3}}\left[ {\sum\limits_{m = 1}^\infty  {\frac{{{a_m}}}
{{\max \{ {m^\lambda },{n^\lambda }\} }}} } \right]} ^p} < k_{_\lambda }^{p - 1}\sum\limits_{n = 1}^\infty  {\left[ {{k_\lambda } - \frac{1}
{{3(p + \lambda  - 2){n^{(p + \lambda  - 2)/p}}}}} \right]{n^{(p - 1)(2 - \lambda ) - 1}}a_n^p}, \label{3.4}
 \tag{3.4}
\end{equation*}
Where $k_\lambda = \frac{\lambda pq}{(p + \lambda - 2)(q + \lambda
- 2)} > 0$. \end{cor}

Taking $\lambda = 1$, in (3.1) and (3.2), we have:

\begin{cor}
\label{Corollary 3.2} If $p,q>1$,$\textstyle{1 \over p} +
\textstyle{1 \over q} = 1$,$a_n \ge 0$, $b_n \ge 0$, such that $0
< \sum\limits_{n = 1}^\infty {n^{p - 2}a_n^p } < \infty $, and $0
< \sum\limits_{n = 1}^\infty {n^{q - 2}b_n^q } < \infty $, then
\begin{equation*}
\begin{split}
&\sum\limits_{n = 1}^\infty  {\sum\limits_{m = 1}^\infty  {\frac{{{a_m}{b_n}}}
{{\max \{ m,n\} }}} }  < pq{\left\{ {\sum\limits_{n = 1}^\infty  {\left[ {1 - \frac{1}
{{3q(p - 1){n^{(p - 1)/p}}}}} \right]{n^{p - 2}}a_n^p} } \right\}^{\tfrac{1}
{p}}}\\
  &\quad \quad \quad \quad \quad \quad \quad \quad  \times {\left\{ {\sum\limits_{n = 1}^\infty  {\left[ {1 - \frac{1}
{{3p(q - 1){n^{(q - 1)/q}}}}} \right]{n^{q - 2}}b_n^q} } \right\}^{\tfrac{1}
{q}}}, \\
 \end{split}
 \label{3.5}\tag{3.5}
\end{equation*}
\begin{equation*}
{\sum\limits_{n = 1}^\infty  {{n^{p - 2}}\left[ {\sum\limits_{m = 1}^\infty  {\frac{{{a_m}}}
{{\max \{ m,n\} }}} } \right]} ^p} < {(pq)^p}\sum\limits_{n = 1}^\infty  {\left[ {1 - \frac{1}
{{3q(p - 1){n^{(p - 1)/p}}}}} \right]{n^{p - 2}}a_n^p}.
 \label{3.6}\tag{3.6}
\end{equation*}
 \end{cor}

In particular, for $p = q = 2$,we have
\begin{equation*}
\sum\limits_{n = 1}^\infty  {\sum\limits_{m = 1}^\infty  {\frac{{{a_m}{b_n}}}
{{\max \{ m,n\} }}} }  < 4{\left\{ {\sum\limits_{n = 1}^\infty  {\left[ {1 - \frac{1}
{{6\sqrt n }}} \right]a_n^2} } \right\}^{\tfrac{1}
{2}}}{\left\{ {\sum\limits_{n = 1}^\infty  {\left[ {1 - \frac{1}
{{6\sqrt n }}} \right]b_n^2} } \right\}^{\tfrac{1}
{2}}},
 \label{3.7}\tag{3.7}
\end{equation*}
\begin{equation*}
{\sum\limits_{n = 1}^\infty  {\left[ {\sum\limits_{m = 1}^\infty  {\frac{{{a_m}}}
{{\max \{ m,n\} }}} } \right]} ^2} < 16\sum\limits_{n = 1}^\infty  {\left[ {1 - \frac{1}
{{6\sqrt n }}} \right]a_n^2}.
 \label{3.8}\tag{3.8}
\end{equation*}


\begin{thebibliography}{99}\small

\bibitem{1}G.H. Hardy, J.E. Littlewood, G. Polya, \textit{Inequalities}, Cambridge Univ.
Press, 1952.

\bibitem{2} D.S. Mitrinovic, J.E. Pecaric, A.M. Fink, \textit{Inequalities Involving
Functions and Their Integrals and Derivatives}, Kluwer Academic
Publishers, Boston, 1991.

\bibitem{3} B.C. Yang, T.M. Rassias, \textit{On the way of weight coefficient and research
for the Hilbert-type inequalities}, Math. Inequal. Appl. 6 (4)
(2003)625--658.

\bibitem{4}B.C. Yang, \textit{Best generalization of Hilbert's type of inequality}, Chinse
Journal of Jilin University (Science edition),2004,42 (1):30-34.

\bibitem{5} B. C.Yang, \textit{The Norm of Operator and Hilbert-type Inequality}, Science
Press, Beijing, 2008.

\bibitem{6} J. Kuang. \textit{Applied Inequalities}, Shandong Science Press, Jinan, 2003.

\end{thebibliography}
 \end{document}